\documentclass{article}
\usepackage{amsthm}
\usepackage{amsfonts}
\usepackage{cite}
\usepackage{amsmath, amscd}
\newtheorem{theorem}{Theorem}[section]
\newtheorem{definition}[theorem]{Definition}
\newtheorem{proposition}[theorem]{Proposition}
\newtheorem{lemma}[theorem]{Lemma}

\newtheorem{corollary}[theorem]{Corollary}

\begin{document}
\author{Jeremy Berquist}
\title{Embeddings of Demi-Normal Varieties}
\maketitle

\noindent
\linebreak
\textbf{Abstract.  }  Our primary result is that a demi-normal quasi-projective variety can be embedded in a demi-normal projective variety.  Recall that a demi-normal variety $X$ is a variety with properties $S_2$, $G_1$, and seminormality.  Equivalently, $X$ has Serre's $S_2$ property and there is an open subvariety $U$ with complement of codimension at least 2 in $X$, such that the only singularities of $U$ are (analytically) double normal crossings.  The term demi-normal was coined by Koll\'ar in \cite{Kol13}.  As a consequence of this embedding theorem, we prove a semi-smooth Grauert-Riemenschneider vanishing theorem for quasi-projective varieties, the projective case having been settled in \cite{Berq14}.  The original form of this vanishing result appears in \cite{GR70}.  We prove an analogous result for semi-rational singularities.  The definition of semi-rationality requires that the choice of a semi-resolution is immaterial.  This also has been established in the projective case in \cite{Berq14}.  The analogous result for quasi-projective varieties is settled here.  Semi-rational surface singularities have also been studied in \cite{vS87}.  

\newpage
\tableofcontents
\addcontentsline{}{}{}
\newpage
\begin{section}{Introduction}
Demi-normal varieties arise in several contexts.  These are varieties $X$ with properties $S_2$, $G_1$, and seminormality (SN).  Equivalently, $X$ has an open subvariety $U$ such that the complement of $U$ has codimension at least two in $X$, and the only singularities in $U$ are (analytically) double normal crossings.  In other words, at singular closed points $x$ in $U$, the completion $\hat{\mathcal{O}_{X,x}}$ is isomorphic to $k[[x_1, \ldots, x_n]]/(x_1x_2).$  

A primary example is that of a simple normal crossing divisor in a smooth variety.  In particular, in performing a strong resolution of singularities on a variety $X$, there is a proper, birational morphism $f: Y \rightarrow X$ such that $Y$ is smooth, $f$ is an isomorphism over the smooth locus of $X$, and the preimage of the singular locus of $X$ is a simple normal crossing divisor in $Y$.  This preimage, viewed as a scheme in its own right, is demi-normal.

The properties $G_1$ and $S_2$ are typically the weakest assumed for non-normal varieties that arise in birational classification.  With these properties, it is possible to develop a theory of Weil divisors on non-normal varieties \cite{Hart94}.  Note that a normal variety is one with properties $R_1$ and $S_2$ and is automatically seminormal.  As $R_1$ is stronger that $G_1$, every normal variety is automatically demi-normal.

At least two results that hold for normal varieties hold also for demi-normal varieties, at least when resolution of singularities is replaced by semiresolution.  One of these is Grauert-Riemenschneider (G-R) vanishing.  If $X$ is a normal variety and $f : Y \rightarrow X$ is a proper, birational morphism from a smooth variety to $X$, then the higher direct images $R^if_*\omega_Y$ are 0 ($i > 0$).  If $X$ is demi-normal and $f$ is a semiresolution (that is, $Y$ is semi-smooth), the same result holds.  This was proved in the case that $X$ is projective in \cite{Berq14}.  We extend this result to the quasi-projective case presently.

A normal variety has rational singularities if for any morphism $f$ as above, the higher direct images $R^if_*\mathcal{O}_Y$ are 0 ($i >0$).  It must be checked that if this is true for one such morphism $f$, then it is true for every such morphism.  Likewise, when $X$ is demi-normal and $f$ is a semiresolution, we say that $X$ has semi-rational singularities when these higher direct images are 0.  Independence of the semiresolution chosen was proved in \cite{Berq14} for projective varieties.  We presently extend the same result to quasi-projective demi-normal varieties.

The key result is that a demi-normal quasi-projective variety can be embedded in a demi-normal projective variety.  Both vanishing results stated above are easier to prove using a projective hypothesis.  In particular, duality and Serre vanishing hold for projective Cohen-Macaulay varieties such as $Y$.  That each of these results holds also for quasi-projective demi-normal varieties requires either a proof that does not depend on such statements for projective varieties, or one that uses the local nature of vanishing together with an embedding result.  We have chosen the latter approach in this paper.

The paper is organized as follows.  In the first section, we recall the definitions of $S_2$, $G_1$, and SN.  We also recall the definition of semiresolution and describe the circumstances under which semiresolutions exist.  

In the second section, we describe the processes of $S_2$-ification and seminormalization.  Both of these end up being ``minimal" in the sense that they preserve the $S_2$ locus and the seminormal locus, respectively.  The goal is to use these processes to replace a projective closure $X$ of our given quasi-projective variety $U$ with a projective closure of $U$ that is both $S_2$ and SN.  We also discuss the construction of the universal pushout.  This is the construction we need to make $X$ Gorenstein in codimension one.

The third section is devoted to the proof of our embedding theorem.  The final section contains corollaries.  In particular, the proof of semi-smooth G-R vanishing and the independence of a chosen semiresolution in the definition of semi-rational singularities appear here.  These results have already been proved in the projective case, and that case is used in the proofs of both results in the quasi-projective case.

All varieties studied here are assumed to be reduced and over an algebraically closed field $k$ of characteristic zero.  In particular, it is possible to speak of resolution of singularities for such varieties.

\end{section}

\begin{section}{Definitions}
We first recall the definitions of $S_2$, $G_1$, and SN.

\begin{definition}  A coherent sheaf $\mathcal{F}$ on a variety $X$ is $S_n$ provided that for all points $x$, we have $$\textnormal{depth } \mathcal{F}_x \geq \textnormal{min}(\textnormal{dim } \mathcal{F}_x, n).$$
\end{definition}

We remark that some authors use a similar definition involving the dimension of the local ring $\mathcal{O}_{X,x}$ in place of the dimension of the module $\mathcal{F}_x$.  The difference between the two definitions is in whether one wants to think of a module as a module over a ring or over the quotient via the annihilator of the module.  As we will call on this alternate definition (without any ambiguity as to which definition is serving which purpose!), for the purposes of this paper let us call it ``strong $S_n$":

\begin{definition}  A coherent sheaf $\mathcal{F}$ on a variety $X$ is strong $S_n$ provided that for all points $x$, we have $$\textnormal{depth } \mathcal{F}_x \geq \textnormal{min}(\textnormal{dim } \mathcal{O}_{X,x}, n).$$
\end{definition}

Of course, here strong $S_n$ implies $S_n$.

A coherent sheaf is \textit{Cohen-Macaulay} if it is $S_n$ for all $n$.  Thus, considering that the depth of a module is always bounded above by its dimension, a Cohen-Macaulay sheaf is one for which the depth and dimension are the same at each localization.

Of course, a variety will be called $S_n$ (or Cohen-Macaulay) if its structure sheaf $\mathcal{O}_X$ has this property.

\begin{definition}  A variety is said to have condition $G_1$ if $\mathcal{O}_{X,x}$ is a Gorenstein ring whenever $x$ is a point of codimension 0 or 1.  In other words, at such points, the canonical module is trivial and the variety is Cohen-Macaulay.
\end{definition} 
\noindent
The canonical module will be discussed in more detail later; we just note that the dualizing sheaf exists for all of our varieties, and that the localizations of the dualizing sheaves are the canonical modules for the corresponding local rings.  Thus, $G_1$ for us simply means that the dualizing sheaf is invertible in codimensions 0 and 1.

We have the following implications for local rings:  $$\textnormal{regular} \implies \textnormal{Gorenstein} \implies \textnormal{Cohen-Macaulay}.$$

\noindent
In fact, any complete intersection is Gorenstein.  In particular, hypersurfaces are Gorenstein, hence $G_1$ and $S_2$.

\begin{definition}  An extension of rings $A \hookrightarrow B$ is a quasi-isomorphism if it is finite, a bijection on prime spectra (hence a homeomorphism on prime spectra), and each residue field extension $k(p) \hookrightarrow k(q)$ is an isomorphism.
\end{definition}
\noindent
The term ``subintegral" is also used in place of ``quasi-isomorphism."

\begin{definition}  Given a finite extension of rings $A \hookrightarrow B$, the seminormalization of $A$ in $B$ is the unique largest subring of $B$ that is quasi-isomorphic to $A$.  We say that $A$ is seminormal in $B$ if it equals its seminormalization in $B$.  We say that a reduced ring $A$ is seminormal if it equals its seminormalization in its integral closure.
\end{definition}

We note that the normalization is finite for reduced finitely-generated algebras over a field.  It is not hard to show that $A$ is seminormal in $B$ if $b^2, b^3 \in A$ imply $b \in A$ for any $b \in B.$

We say that a (reduced) variety is seminormal if each of its local rings is seminormal.

\begin{proposition}[\cite{Kol13}, 5.1] A variety $X$ is $S_2$, $G_1$, and SN if and only if it is $S_2$ and there exists an open subvariety $U$ such that:  (i)  the complement of $U$ in $X$ has codimension at least two; (ii)  for any closed singular point $x \in U$, $\hat{\mathcal{O}}_{X,x} \cong k[[x_1, x_2, \ldots, x_n]]/(x_1x_2).$  
\end{proposition}

We now come to the definition of a semiresolution.  Recall that a pinch point is a point whose local ring is analytically isomorphic to $k[[x_1, x_2, \ldots, x_n]]/(x_1^2-x_2^2x_3).$  The pinch point is significant in that blowing up the origin produces another pinch point, and thus cannot be simplified without blowing up the entire double locus.  A pinch point on a surface is a quotient of a double normal crossing point.  See [7, 10.4] for more details on the relationship between double normal crossings and pinch points.  In the surface case, both are examples of the same phenomenon.  See also \cite{vS87}.

\begin{definition}
A variety is semismooth if every closed point is either smooth, a double normal crossing point, or a pinch point.
\end{definition}

\begin{definition}  Let $X$ be a variety as above.  Then a morphism $f: Y \rightarrow X$ is a semiresolution if the following conditions are satisfied:  (i)  $f$ is projective, (ii) $Y$ is semismooth, (iii) no component of the conductor $C_Y$ is $f$-exceptional, and (iv) $f|_{f^{-1}(U)}: f^{-1}(U) \rightarrow U$ is an isomorphism, where $U$ is an open set whose only closed singular points are double normal crossings.
\end{definition}

If $U \subseteq X$ is an open subvariety whose closed singular points are (analytically) double normal crossing points or pinch points, then Koll\'ar shows that there is a projective morphism $f: Y \rightarrow X$ such that $Y$ is semismooth, $f^{-1}(U) \rightarrow U$ is an isomorphism, and the singular (double) locus of $Y$ maps birationally onto the singular locus of $U$.  Thus semiresolutions exist.  We should stress that the characteristic of the base field is zero; Koll\'ar's construction uses a resolution of singularities of the normalization of $X$.

\end{section}
\begin{section}{Constructions}
We now discuss the constructions used in proving the main embedding theorem.  

Recall first that a reduced ring $A$ has a seminormalization in its integral closure $\overline{A}$.  See 2.5.  Seminormalization is a local procedure.  In particular, when $B$ is the seminormalization of $A$, $$A \hookrightarrow B \hookrightarrow \overline{A},$$ then at every localization $A_p$, $B_p$ is the seminormalization of $A_p$.  See \cite{GT80}, (2.1).  Moreover, when $A_p$ is seminormal, already $A_p = B_p$.  

The good part is that we can patch together local data to obtain a seminormalization of a variety, which is an isomorphism over the seminormal locus.  As long as $A$ is of finite type over a field $k$, the normalization $\overline{A}$ is a finite $A$-module, hence so is any intermediate ring.  That is, the seminormalization is finite as an $A$-module.  Specifically, we have:  

\begin{proposition}  Let $X$ be a reduced variety.  Then there is a finite, birational morphism $f: X_{SN} \rightarrow X$ such that $X_{SN}$ is seminormal, and such that $f$ is an isomorphism over the seminormal locus of $X$.  Moreover, $f$ is a bijection on points, and all induced residue field extensions are isomorphisms.
\proof  This proceeds along the same lines as \cite{Hart87}, II.Ex.3.8.  We omit a full proof.  That $f$ is a bijection on points, hence a homeomorphism, and that all residue field extensions are isomorphisms follows from the fact that both are characteristics of the seminormalization of a reduced ring.  \qed
\end{proposition}

Our second construction is that of $S_2$-ification.  As the name suggests, we wish to construct a finite ring extension $A \hookrightarrow B$ such that $B$ is $S_2$ as a ring.  Here $A$ is reduced and of finite type over a field $k$.  For this construction, we will need to use the canonical module $\omega_A$.

One way to define $\omega_A$ is via noether normalization.  When $A$ is of finite type over a field $k$, then there is a polynomial ring $C$ over $k$ contained in $A$ such that $A$ is finite as a $C$-module.  Then duality for a finite morphism implies that $\omega_A = \textnormal{Hom}_C(A, \omega_C).$  A polynomial ring has localizations that are regular local rings.  Thus $\omega_C$ is isomorphic to $C$ as a $C$-module.  So we have $\omega_A = \textnormal{Hom}_C(A, C)$.  This implies, among other things, that $\omega_A$ is finite and strong $S_2$ as an $A$-module, as the next two lemmas show.

\begin{lemma}  If $M$ and $N$ are finite $A$-modules, then so is $\textnormal{Hom}_A(M,N)$.  If moreover $N$ is strong $S_2$ as an $A$-module, then so is $\textnormal{Hom}_A(M,N)$.
\proof  Write $M = Au_1 + \cdots + Au_m$ and $N = Av_1 + \cdots + Av_n$.  Then define $\phi_{i,j}: M \rightarrow N$ by $\phi_{i,j}(u_i) = v_j$ and $\phi_{i,j}(u_k) = 0$ for $i \neq k$.  For an arbitrary $\phi: M \rightarrow N$, write $\phi(u_i)= \Sigma a_{ij}v_j$.  Then $\phi = \Sigma a_{ij}\phi_{i,j}$.  Thus $\textnormal{Hom}_A(M,N)$ is spanned by the $\phi_{i,j}$.

Suppose next that $N$ is strong $S_2$ as an $A$-module.  Suppose that $A_p$ has dimension at least two.  We prove that $\textnormal{Hom}_A(M,N)$ has depth at least two at $p$.  The rest of the proof that the Hom module is strong $S_2$ is proved similarly; that is, it would remain to check the cases in which $A_p$ has dimension zero (which is trivial) or one (which proof is included in what follows).  

We know that $N_p$ has depth at least two.  Since $M$ is finite as an $A$-module, the localization of the Hom module at $p$ is isomorphic to $\textnormal{Hom}_{A_p}(M_p, N_p)$.  Suppose that $z_1, z_2$ is an $N_p$-regular sequence.  We show that these elements are also regular for the Hom module at $p$.

Suppose first that $z_1\phi = 0$ for some $\phi \in \textnormal{Hom}_{A_p}(M_p, N_p).$  Then $z_1\phi(u) = 0$ for all $u \in M_p$.  Since $z_1$ is $N_p$-regular, this implies that $\phi(u) = 0$ for all $u$, meaning that $\phi$ is identically zero.  So $z_1$ is $\textnormal{Hom}_{A_p}(M_p,N_p)$-regular.

Next, suppose that $z_2\phi = z_1\psi$ for morphisms $\phi$ and $\psi$.  Evaluating at an arbitrary element of $M_p$, we again see that $\phi(u) = z_1v$ for some $v \in N_p$.  This element $v$ depends uniquely on $u$ since $z_1$ is $N_p$-regular.  Define $\rho(u) = v$.  Then $\rho$ is an $A_p$-module homomorphism, and by construction $\phi = z_1\rho$.  This completes the proof that $z_1, z_2$ is a regular sequence for the Hom module, hence that the depth is at least two.\qed
\end{lemma}

What we have shown is that $\omega_A$ is finite and strong $S_2$ as a $C$-module.  Since $A$ is finite over $C$, $\omega_A$ is also finite as an $A$-module.  That it is strong $S_2$ as an $A$-module follows from the next lemma.

\begin{lemma}  Let $C \hookrightarrow A$ be a finite ring extension and let $M$ be an $A$-module that is strong $S_2$ as a $C$-module.  Then $M$ is strong $S_2$ as an $A$-module.
\proof  For simplicity, we again consider just the case in which $A_p$ has dimension at least two.  We wish to show that $M_p$ has depth at least two.  

Let $q = C \cap p$.  Then obviously the height of $q$ is at least two.  Thus $M_q$, considered as a $C_q$-module, has depth at least two.  There are elements $z_1, z_2 \in qC_q$ that form an $M_q$-regular sequence.  In particular, $\dim(M_q/(z_1, z_2)) < \dim(M_q/(z_1)) < \dim(M_q)$ as $C_q$-modules.  

Clearly the elements $z_1, z_2$ are in $pA_p$. Moreover, $M_p$ is a localization of $M_q$, and localizing preserves the dimension since $p$ sits over $q$.  Thus $$\dim(M_p/(z_1, z_2)) < \dim(M_p/(z_1)) < \dim(M_p)$$ as $A_p$-modules.  The dimension can go down only for regular elements.  Thus $M_p$ also has depth at least two, as we wanted to show.  \qed
\end{lemma}

In light of these two lemmas, we consider the $A$-module $\textnormal{Hom}_A(\omega_A, \omega_A).$  Since $\omega_A$ is strong $S_2$ as an $A$-module, so is the Hom module.  Moreover, $\textnormal{Hom}_A(\omega_A, \omega_A)$ is a ring whose multiplicative structure is given by function composition.  It is therefore $S_2$ as a ring and finite as an $A$-module.  Since $\omega_A$ is torsion-free (by its definition), we conclude that $$A \hookrightarrow \textnormal{Hom}_A(\omega_A, \omega_A)$$ represents an $S_2$-ification of $A$.

We would like this to have a ``minimality" property.  That is, we would like for the Hom module to embed in $B$ for any finite extension $A \hookrightarrow B$ in which $B$ is $S_2$.  This would imply, in particular, that $A$ is $S_2$ if and only if the above inclusion is an equality.

This is, in fact, true, but not obvious.  It turns out that $\omega_A$ is self-dualizing on the set of strong $S_2$ $A$-modules.  In other words, when $B$ is strong $S_2$ as an $A$-module, the natural inclusion $$B \hookrightarrow \textnormal{Hom}_A(\textnormal{Hom}_A(B, \omega_A), \omega_A)$$ is an equality.

Now consider that $A \hookrightarrow B$ is a finite extension, where $B$ is $S_2$.  Then $A$ maps into $$\textnormal{Hom}_A(\omega_A, \omega_A) = \textnormal{Hom}_A(\textnormal{Hom}_A(A, \omega_A), \omega_A),$$ which injects naturally into the double dual of $B$.  This is what we wanted to show. 

Finally, we observe that any two $S_2$-ifications of $A$ are naturally isomorphic.  This is just the minimality property discussed above.  It follows that we can patch together $S_2$-ifications for local data corresponding to an open cover of a given variety $X$.  In particular, we can state the following:

\begin{proposition}  Let $X$ be a reduced variety.  Then there is a finite, birational morphism $f: X_{S_2} \rightarrow X$ such that $X_{S_2}$ is $S_2$ and such that $f$ is an isomorphism over the $S_2$ locus of $X$.
\proof  This is the content of the above discussion.\qed
\end{proposition}

It is our goal to use both seminormalization and $S_2$-ification in order to produce a variety with both properties.  For this, it is essential that there is the following preservation property.

\begin{lemma}  Let $A$ be a seminormal ring, and let $C$ be a subring of $\overline{A}$ containing $A$.  If $C$ is $S_2$, then it is seminormal.
\proof  This is \cite{GT80}, (2.8).\qed
\end{lemma}

Therefore, starting with our ring $A$, we perform seminormalization, obtaining $A \hookrightarrow B$, where $B$ is finite as an $A$-module and contained in $\overline{A}$.  Then we perform $S_2$-ification, obtaining $A \hookrightarrow B \hookrightarrow C$, where $C$ is $S_2$ and finite as a $B$-module and contained in $\overline{A}$.  We compose to obtain a finite $A$-module contained in $\overline{A}$, which is both seminormal and $S_2$. 

These are local and birational operations.  We can combine (3.1) and (3.4) to obtain a finite, birational morphism $f: Y \rightarrow X$ such that $Y$ is seminormal and $S_2$, and such that $f$ is an isomorphism over the locus in $X$ which is both seminormal and $S_2$.

The next section uses these results in the proof of the main embedding theorem.  We will need one final lemma, which states that the $Y$ obtained above is projective when $X$ is projective.

\begin{lemma}  Let $f: Y \rightarrow X$ be finite, where $X$ is projective.  Then $Y$ is projective.
\proof  Both finite and projective morphisms are proper, as is a composition of proper morphisms.  Thus $Y$ is proper.  It is projective, provided that there exists an ample sheaf on $Y$.  We know that there is a (very) ample sheaf $\mathcal{L}$ on $X$.  We prove that $f^*\mathcal{L}$ is ample on $Y$.

Let $\mathcal{F}$ be coherent on $Y$.  Then $f_*\mathcal{F}$ is coherent on $X$, so $f_*\mathcal{F} \otimes \mathcal{L}^n$ is globally generated for some $n \gg 0$.  There is a surjective morphism $\oplus_{i \in I} \mathcal{O}_X \rightarrow f_*\mathcal{F} \otimes \mathcal{L}^n \rightarrow 0$.  Now $f^*$ is right exact and preserves direct sums, since the same is true for tensor product, and $f$ is finite.  Applying $f^*$ to this surjection, we obtain $\oplus_{i \in I} \mathcal{O}_Y \rightarrow f^*f_*\mathcal{F} \otimes f^*\mathcal{L}^n \rightarrow 0$.  Since $f^*\mathcal{L}^n$ is invertible, we will be done if we can show that there is a surjection $f^*f_*\mathcal{F} \rightarrow \mathcal{F} \rightarrow 0$.

Since $f$ is finite, this reduces to the fact that there is a surjective morphism $M \otimes_A B \rightarrow M$, where $M$ is a finite $B$-module and $A \rightarrow B$ is a ring extension.  This surjection is simply given by multiplication $m \otimes b \mapsto b \cdot m$.  This finishes the proof.\qed
\end{lemma}
\end{section}

\begin{section}{Proof of the Embedding Theorem}
We now come to the main result.

\begin{theorem}  Let $U$ be a demi-normal quasi-projective variety.  Then $U$ can be embedded as a dense open set in a demi-normal projective variety.
\proof  Let $U \hookrightarrow X$ be a projective closure of $U$.  By possibly replacing $X$ with the closure of $U$, we can assume that $U$ is dense in $X$.  Let $g:  X_1 \rightarrow X$ be the normalization of $X$.  Note that $X_1$ is projective, by 3.6.  Let $C$ be the conductor in $U$, and let $\overline{C}$ be its closure in $X$.  Note that the only codimension one points in $\overline{C}$ are already in $C$; that is, if $p$ is a codimension one point of $X$ outside of $U$, then $\overline{C} \cap \{ p\}^-$ has codimension at least two in $X$.

 According to \cite{Art70} 3.1, the universal pushout can be used to glue along a morphism without affecting the rest of a given variety.  Specifically, given a closed subscheme $B \hookrightarrow Y$ and a finite morphism $B \rightarrow B/\tau$, the universal pushout $$\begin{CD}
B           @>>>   Y  \\        
@VVV               @VVV  \\         
B/\tau      @>>>   Y' 
\end{CD}$$
has the property that $Y \rightarrow Y'$ is proper, agrees with $B \rightarrow B/\tau$ on $B$ and is an isomorphism elsewhere.

We consider the universal pushout $$\begin{CD}
g^{-1}(\overline{C}) @>>> X_1 \\
@VgVV                      @VpVV \\
\overline{C} @>>>         X_2 \\
\end{CD}$$

\noindent
Here $p$ is proper and birational, $X_2$ contains a copy of $\overline{C}$, and $X_1 - g^{-1}(\overline{C}) \rightarrow X_2 - \overline{C}$ is an isomorphism.  

Observe that $C = U \cap \overline{C}$.  In the universal pushout $$\begin{CD}
g^{-1}(C) @>>> g^{-1}(U) \\
@VVV           @VVV \\
C         @>>> U \\
\end{CD}$$

\noindent
$g^{-1}(U) \rightarrow U$ is the normalization of $U$, hence $U$ is obtained in the usual way by gluing along the conductor.  Moreover, since $g^{-1}(C) = g^{-1}(U) \cap g^{-1}(\overline{C})$ is open in $g^{-1}(\overline{C})$, $U$ embeds as an open subvariety of $X_2$.  Now all the codimension one points of $\overline{C}$ in $X_2$ are in $U$, and outside $\overline{C}$, $X_2$ looks like $X_1$.  Thus $X_2$ is regular in codimension one outside of $U$, and therefore is $G_1$.

We next observe that $X_2$ is projective.  In fact, there is a finite morphism from $X_2$ to $X$, and we can apply 3.6.  To show that there is such a morphism, we observe first that $X$ is obtained by gluing along the conductor in its normalization $X_1$.  Moreover, if $D$ is the conductor in $X$, then $\overline{C}$ is contained in $D$.  In other words, $D^c$ is the set of normal points in $X$, so that $D^c \cap U$ consists of normal points of $U$, so that $D^c \cap U \subset C^c$.  Thus $C \subset D$, and since $D$ is closed, $\overline{C} \subset D$.  

By the universal property of the pushout, there is an induced morphism $X_2 \rightarrow X$.  To see that this morphism is finite, it is enough to look locally.  Suppose that $X$ is given locally by Spec$ A$.  Then on the ring level, $X_2$ is given by Spec$ B$, where $B$ is a universal pullback $$\begin{CD}
B @>>> \overline{A} \\
@VVV   @VVV \\
A/J @>>> \overline{A}/J\overline{A} \\
\end{CD}$$

\noindent
where $\overline{A}$ is the normalization of $A$.  So $A$ maps naturally into $B$, and we have inclusions $A \hookrightarrow B \hookrightarrow \overline{A}$.  Since $\overline{A}$ is finite as an $A$-module, so is $B$.  So $X_2$ is in fact finite over $X$, hence projective.  Note again that $X_2$ is birational to $X$, since $X_2$ is birational to the normalization $X_1$ of $X$.  In particular, $U$ is dense in $X_2$.

We now have an inclusion $U \hookrightarrow X_2$, where $X_2$ is $G_1$.  Using 3.1, 3.4, and 3.5, there is a finite, birational morphism $f:  X_3 \rightarrow X_2$, such that $X_3$ is seminormal and $S_2$, and such that $f$ is an isomorphism over $U$.  In fact, $f$ is an isomorphism over all codimension one points of $X_2$, since outside of $U$, $X_2$ is regular in codimension one outside of $U$, and regular implies both seminormal (even normal) and $S_2$ (even Cohen-Macaulay).  

In conclusion, $X_3$ is seminormal, $S_2$, and $G_1$.  Hence it is demi-normal.  Since $X$ is projective and all morphisms are finite, $X_3$ is also projective.  Finally, since all morphisms are birational, $U$ embeds as a dense open subset of $X_3$.  This completes the proof.\qed
\end{theorem}

\noindent
\textbf{Remark.  }  Note that the proof of 4.1 shows also that a normal, quasi-projective variety can be embedded in a normal, projective variety.  In fact, normality is determined by the conditions $R_1$ and $S_2$.  If $U$ has these properties, then the proof of 4.1 shows that $U$ can be embedded in a projective variety that is regular in codimension one outside of $U$ and $S_2$.

\noindent
\linebreak
In the next section, we give a few corollaries of the embedding result.  In particular, we improve two results from \cite{Berq14}.  Namely, we show a G-R vanishing theorem for quasi-projective varieties and an independence result concerning semi-rational singularities, also for quasi-projective varieties, the projective cases having been settled in \cite{Berq14}.
\end{section}
\begin{section}{Corollaries}
Our first corollary of 4.1 is a semi-smooth Grauert-Riemenschneider vanishing theorem.

\begin{corollary}  Let $U$ be a quasi-projective, demi-normal variety.  Let $f: V \rightarrow U$ be a semiresolution.  Then $R^if_*\omega_V = 0$ for $i > 0$.
\proof  We use the fact that this same result holds when $U$ is replaced by a projective demi-normal variety, as shown in \cite{Berq14}.  Let $U \hookrightarrow X$ be a projective closure such that $X$ is demi-normal.  Noting that $V$ is also quasi-projective (see \cite{Hart87}, II.7.16 and II.7.17), let $V \hookrightarrow Y$ be a projective closure such that $Y$ is demi-normal.  Then there is a birational map $Y \dashrightarrow X$.  The domain of definition includes $V$.  

Now consider the diagram $$\begin{CD}
V      @>>>   Y \\
@VVV          @VVV \\
U      @>>>   X \\
\end{CD}$$

The right arrow represents only a birational map. Since $V$ is in the domain of definition, when we resolve the indeterminacies, which amounts to blowing up the complement of the domain of definition, $V$ will still embed in the resulting variety.  So let $Z_0 \rightarrow Y$ be a sequence of blowups so that $Z_0 \rightarrow X$ is a morphism.  We obtain a diagram as before, with $Z_0$ in place of $Y$, where now the right map is a morphism.  

Here $Z_0$ is a projective closure of $V$, but is not necessarily demi-normal.  There is a finite, birational morphism $Z_1 \rightarrow Z_0$, which is an isomorphism over $V$, and such that $Z_1$ is demi-normal.  That is the content of the proof of 4.1.

Now we let $Z_2 \rightarrow Z_1$ be a semiresolution, which we can assume to be an isomorphism over $V$ since $V$ is semismooth.  Then $Z_2$ is projective, and we obtain a diagram as before with $Z_2$ in place of $Y$.  In other words, we have a commutative diagram $$\begin{CD}
V @>>> Z_2 \\
@VVV   @VVV \\
U @>>> X \\
\end{CD}$$

\noindent 
where $g:  Z_2 \rightarrow X$ is a projective, birational morphism from a semismooth variety to $X$. By the G-R vanishing theorem for projective varieties, $R^ig_*\omega_{Z_2} = 0$ for $i >0$.  Now since $V \hookrightarrow Z_2$ is an open immersion, $\omega_V = \omega_{Z_2}|_V.$  We obtain the desired result by observing that $R^if_*\omega_V \cong R^ig_*\omega_{Z_2}|_U.$\qed
\end{corollary}

As a second corollary, we show a result used in showing that semi-rational singularities are well-defined.  The key observation in that result for projective varieties is that higher direct images are zero for a morphism between semismooth projective varieties.

\begin{corollary}  Let $f: V \rightarrow U$ be a projective, birational morphism between quasi-projective, semismooth varieties.  Then $R^if_*\mathcal{O}_V = 0$ for $i > 0$.
\proof  The proof proceeds along the same lines as 5.1.  The important part is that a quasi-projective, semismooth variety $U$ embeds into a projective semismooth variety.  This was shown already.  If $U \hookrightarrow X$ is an embedding such that $X$ is demi-normal, then there is a semiresolution $Y \rightarrow X$ that is an isomorphism over $U$, since $U$ is semismooth.  If we embed both $U$ and $V$ this way, then we are in the situation of 5.1.  The vanishing of the higher direct images of $\mathcal{O}_V$ follows from the same arguments used to show that the higher direct images of $\omega_V$ vanish for a morphism such as $f$.\qed
\end{corollary}

\noindent
\textbf{Remark.  }  It seems natural to ask whether a quasi-projective variety with (semi-) rational singularities can be embedded in a projective variety with (semi)-rational singularities.  This is an open problem.  In fact, it is not known even in the case where $X$ is 3-dimensional and has only one non-rational point.  Suppose $X$ is singular along a curve $C$ and has rational singularities along $C\ \{p\}$ for a single point $p$.  It is not known whether one can partially resolve $p$ to become a rational singularity without resolving the singularities along $C$.

\end{section}
\nocite{Art70}
\nocite{Kol13}
\nocite{Hart94}
\nocite{Hart87}
\nocite{BH98}
\nocite{GT80}
\nocite{Kov99}
\nocite{LV81}
\nocite{Reid94}
\nocite{dFH09}
\nocite{KSS09}
\nocite{Sch09}
\nocite{vS87}
\nocite{Berq14}
\nocite{GR70}

\bibliographystyle{plain}
\bibliography{paper}

\end{document}